\documentclass{article}

\usepackage{graphicx,color,amsmath,amsfonts,amssymb,times,url}

\begin{document}

\title{Computer geometry: \\ Rep-tiles with a hole}
\markright{Rep-tiles with a hole}
\author{Christoph Bandt and Dmitry Mekhontsev}

\maketitle

\begin{abstract}
A cube is an 8-rep-tile: it is the union of eight smaller copies of itself. Is there a set with a hole which has this property?
The computer found an interesting and complicated solution, which then could be simplified. We discuss some problems of computer-assisted research in geometry. 
\end{abstract}

Will computers help us do geometrical research? Can they find something new? How can we direct them to do those things which we are interested in? On the other hand, will computers change our attitudes?
We discuss such issues for elementary problems of fractal geometry \cite{pinwheel,sierel}, using the free software package IFStile \cite{M}. This note is concerned with self-similar tilings of three-dimensional space.

\begin{figure}[h]
	\centerline{
		\includegraphics[height=0.28\textwidth]{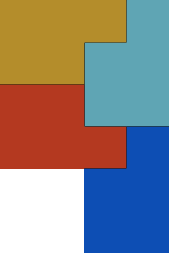}\quad\
		\includegraphics[height=0.3\textwidth]{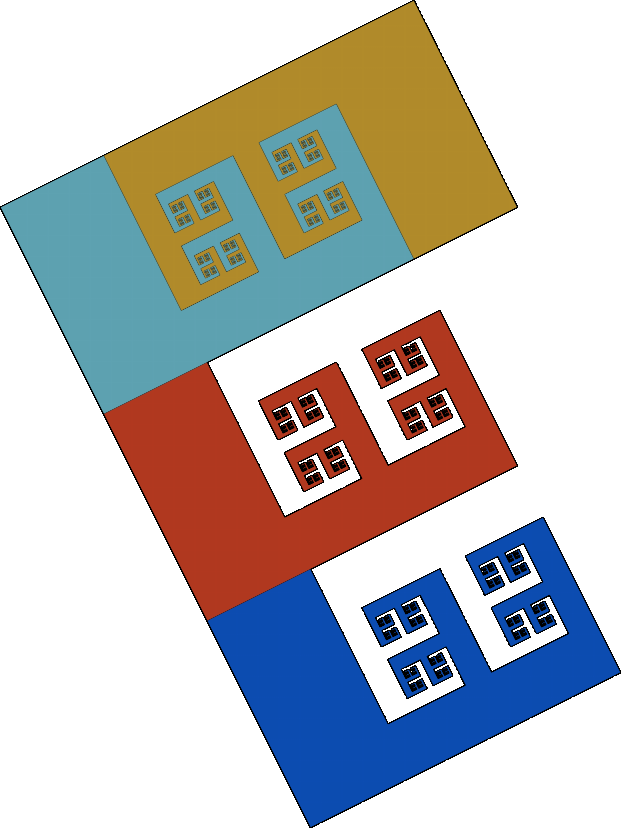}\quad\ 
             \includegraphics[height=0.3\textwidth]{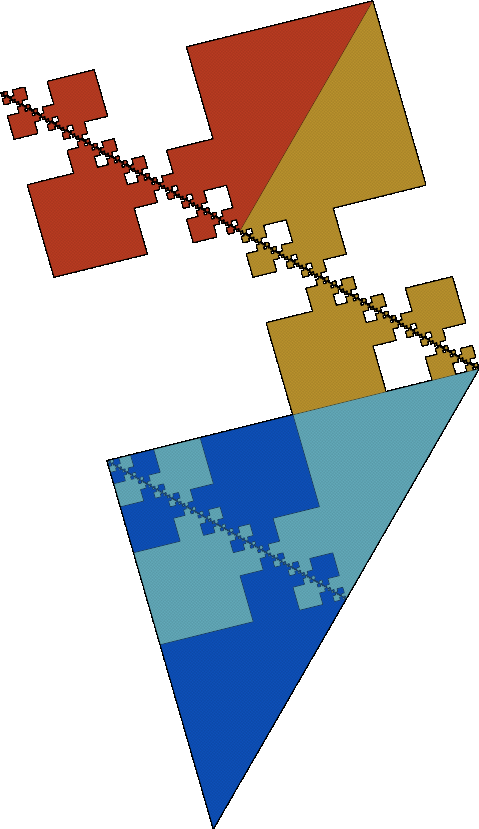}\quad\
             \includegraphics[height=0.3\textwidth]{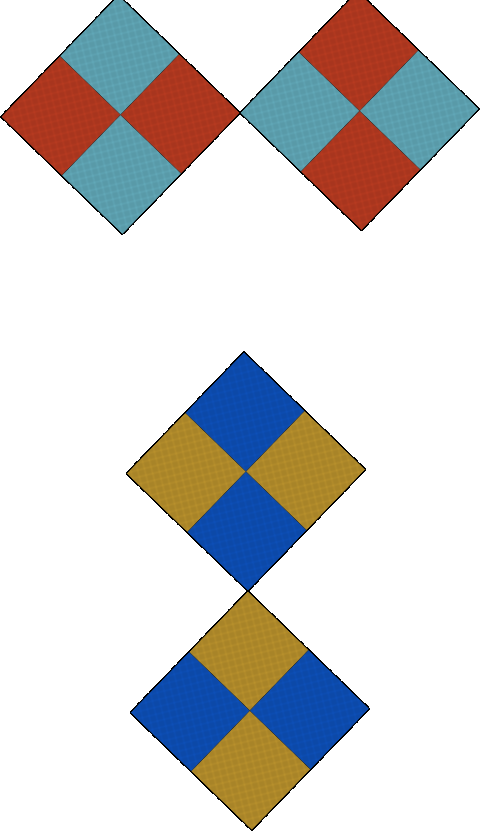}
	}
	\caption{4-rep-tiles in the plane.}
	\label{rep4}
\end{figure}

\subsection*{Rep-tiles.}
A closed bounded set $A$ with non-empty interior in plane or space is called an $m$-rep-tile if there are sets $A_1, A_2, ..., A_m$ congruent to $A,$ such that different sets $A_k, A_j$ have no common interior points, and the union $B=A_1\cup ...\cup A_m$ is geometrically similar to $A.$ The standard example in the plane is a square, or a parallelogram, or a triangle, with $m=4.$ Some other examples are shown in Figure \ref{rep4}. Exercise: show that a triangle with angles of 30, 60, and 90 degrees is a 3-rep-tile. 

\begin{figure}[h] 
	\centerline{\includegraphics[height=0.43\textwidth]{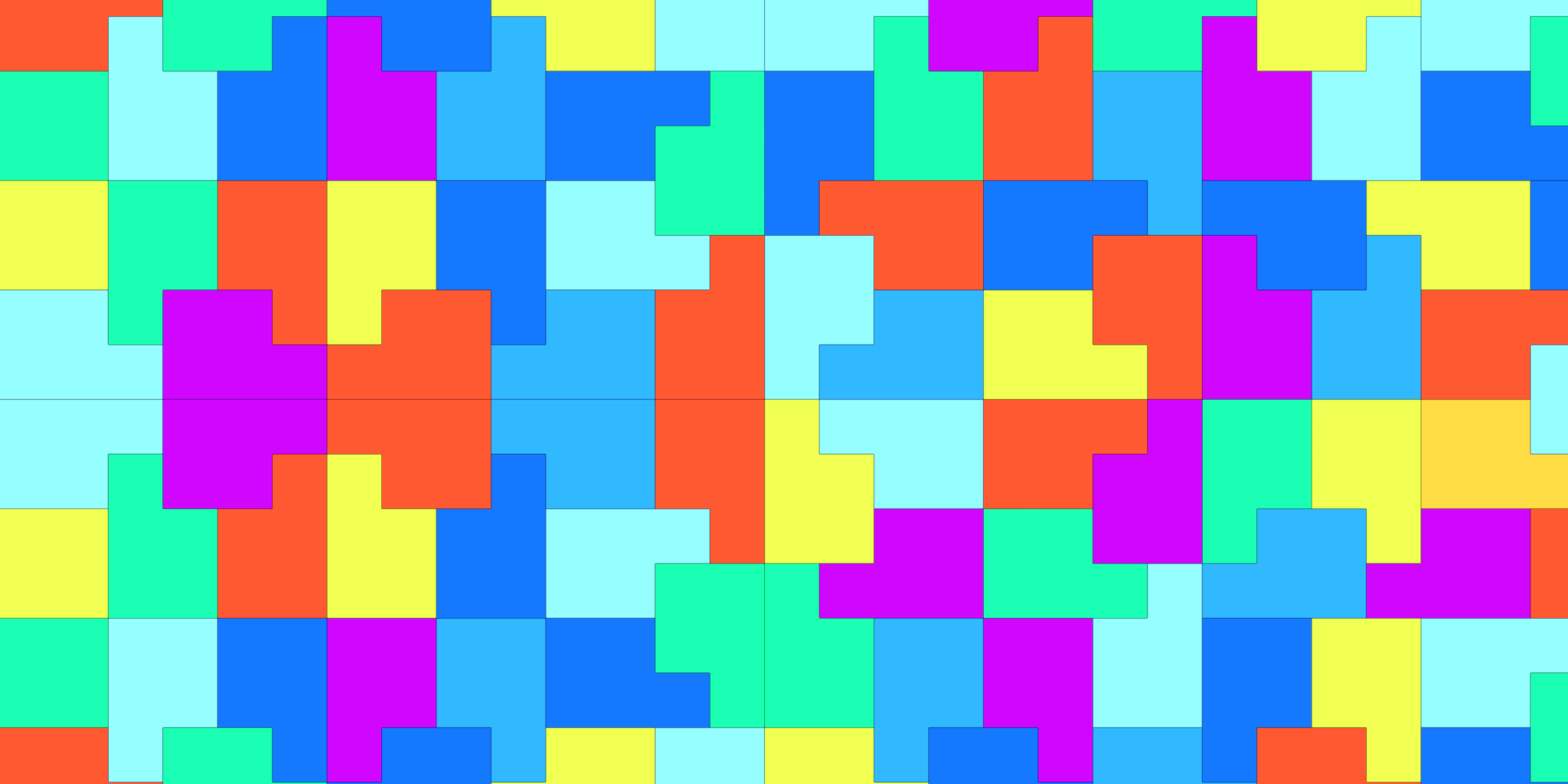}
	}
	\caption{Tiling obtained from the `flag' rep-tile.}
	\label{flagt}
\end{figure}

`Rep' stands for `replication', and the sets are called tiles since they can tile the whole plane. Such tilings can be obtained by observing that $B$ is also a rep-tile and contained in still larger super-rep-tiles $C,D,...,$ and they all are unions of copies of $A$ \cite{GS,Se}. One possible tiling for the `flag' in Figure  \ref{rep4} is indicated in  Figure \ref{flagt}. 
Exercise: try to assemble the tiles in this picture to supertiles. 
The tilings generated by the flag are non-periodic and quite intricate while the $2\times 2$ subdivision of the square provides only the ordinary periodic checkerboard tiling.

\begin{figure}[h] 
	\centerline{
		\includegraphics[height=0.27\textwidth]{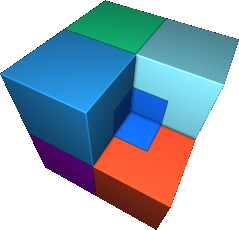}\quad
             \includegraphics[height=0.23\textwidth]{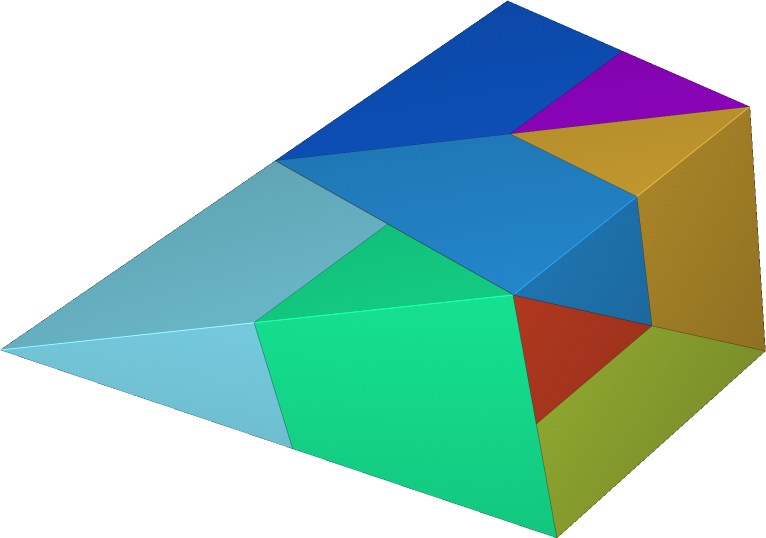}\quad
		\includegraphics[height=0.3\textwidth]{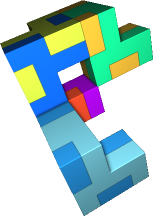}
	}
	\caption{Notched cube and two other polyhedral 8-rep-tiles.}
	\label{rep8}
\end{figure}

Rep-tiles were introduced in 1963 as recreational objects by Gardner \cite{Ga} and Golomb \cite{Gol}. In the 1980s they became interesting as models of quasicrystals  \cite[chapter 11]{GS}, \cite{Se}, as examples of self-similar fractals \cite{Bar}, as a tool for constructing multidimensional wavelets \cite{GM}, and as unit intervals for exotic number systems \cite{Vi}.
For the plane, plenty of $m$-rep-tiles are known for every $m$  \cite{Ba5}. In three-dimensional space, a tetrahedral $m$-rep-tile can exist only for cubic numbers $m,$ not for $m<8$ \cite{MS11}. For $m=8,$ the cube is a standard rep-tile, and the notched cube (`chair') in Figure \ref{rep8} is another well-known example. The regular tetrahedron is not an $8$-rep-tile but some other special tetrahedra are, one of them found by M.J.M. Hill already in 1895, and two others found in 1994 \cite{LJ}. Recent results support the conjecture that there are no further 8-rep-tile tetrahedra \cite{Ha}.  
Figure \ref{rep8} shows two other polyhedral examples found with the IFStile package.

\subsection*{Algebra and algorithms.}
For computer work, geometric concepts must be reformulated in terms of algebra. This was done by John Hutchinson. Instead of saying that $A_k$ is congruent to $A,$ he introduced an isometry map $h_k$ with $A_k=h_k(A).$ Instead of saying that the union $B$ of the copies $A_k$ is geometrically similar to $A,$ he took a similarity mapping $g$ with $B=g(A).$ Of course $g$ must be expanding - it must increase all distances by a factor greater than 1. The defining equation for an $m$-rep-tile becomes
\begin{equation}\label{hut}
g(A)=h_1(A) \cup h_2(A) ... \cup h_m(A) \, ,
\end{equation}
with given data (`coefficients') $g,h_1,...,h_m$ and the unknown set $A.$ Hutchinson proved that this equation always has a unique solution $A$ in the space of compact non-empty subsets of plane or space. His paper \cite{Hu} has become famous just for this rather simple observation although it contains much more difficult theorems. The proof can be found in every textbook on fractal geometry, for example \cite{Bar}.

Since $A$ should have non-empty interior and thus positive volume, a comparison of the volume on both sides of \eqref{hut} shows that $g$ must have determinant $\pm m,$ by a basic theorem of linear algebra. We shall not need this general fact since we consider only the mapping 
\begin{equation}\label{g}
 g(x)=2x\, .
\end{equation}
For this map $g(A)$ has four times larger area than $A$ in the plane, and eight times larger volume in three-dimensional space. So we shall study $m$-rep-tiles with $m=4$ in dimension $d=2$ and with $m=8$ for $d=3.$  Moreover, we consider only isometries $h$ with integer coefficients:
\begin{equation}\label{hk}
 h(x)=M x +v\, ,
\end{equation}
where $v$ is a vector with integer coordinates, and $M$ is a quadratic matrix which has exactly one entry $+1$ or $-1$ in each row and each column, and all other entries are zero. Exercise: there are 8 such matrices for $d=2$ and 48 for $d=3.$
These are the isometries which transform the lattice $\mathbb Z^d$ of integer vectors into itself. The linear maps $f(x)=Mx$ are rotations and reflections which transform the unit cube $[-1,+1]^d$ into itself. These maps form the symmetry group of the unit square for $d=2,$  and of the unit cube for $d=3.$  The resulting rep-tiles form the `square family' and the `cube family', respectively. As we shall see, both families are very large.

It is crucial that our data $g,h_1,...,h_m$ are given by integers!  Integer calculations in the computer are accurate while calculations with real numbers are only approximate, with a numerical error. Extensive calculations are needed since the definition of rep-tile requires that different pieces $A_k, A_j$ of $g(A)$ have no common interior points. To check this condition, we have to study all \textit{neighbor types} $A_k\cap A_j,$ which include also `pieces of pieces' on several levels. They are characterized algebraically by \textit{neighbor maps} which are isometries like the $h_k.$   Exercise: analyse the type of maps (translation, reflection, rotation) which transform tiles in Figure \ref{flagt} into their neighboring tiles. Altogether, there are 60 such maps while for a square tiling we have only 8 translations.

For the case of integer data, there are only finitely many neighbor maps. They can be determined recursively, and if the map $f(x)=x$ is not among them then we really have a rep-tile. Details are explained in \cite{pinwheel, sierel} and the literature quoted there. 
A neighbor map algorithm was implemented in the IFStile package.  It does not only check the rep-tile property, but also calculates the number of neighbor types, as well as the fractal dimension of the boundary, and further parameters which characterize the tile.

Now a search for rep-tiles can be done by randomly generating various data $M_k,v_k$ for $k=1,...,m$ and checking each time whether we obtain a rep-tile. The data and parameters of all resulting rep-tiles will be stored. Within one hour, we get at least 25000 examples. Exercise: download IFStile and try yourself. (Take the square family by clicking the star icon and the first item in the list. For search, click the binocular icon and `Start'.)

Working out a single three-dimensional example by hand may take a day, or even a week. The computer opens up new perspectives. We can explore territories which previously were totally inaccessible to us. 

\begin{figure}[h] 
	\centerline{
		\includegraphics[height=0.3\textwidth]{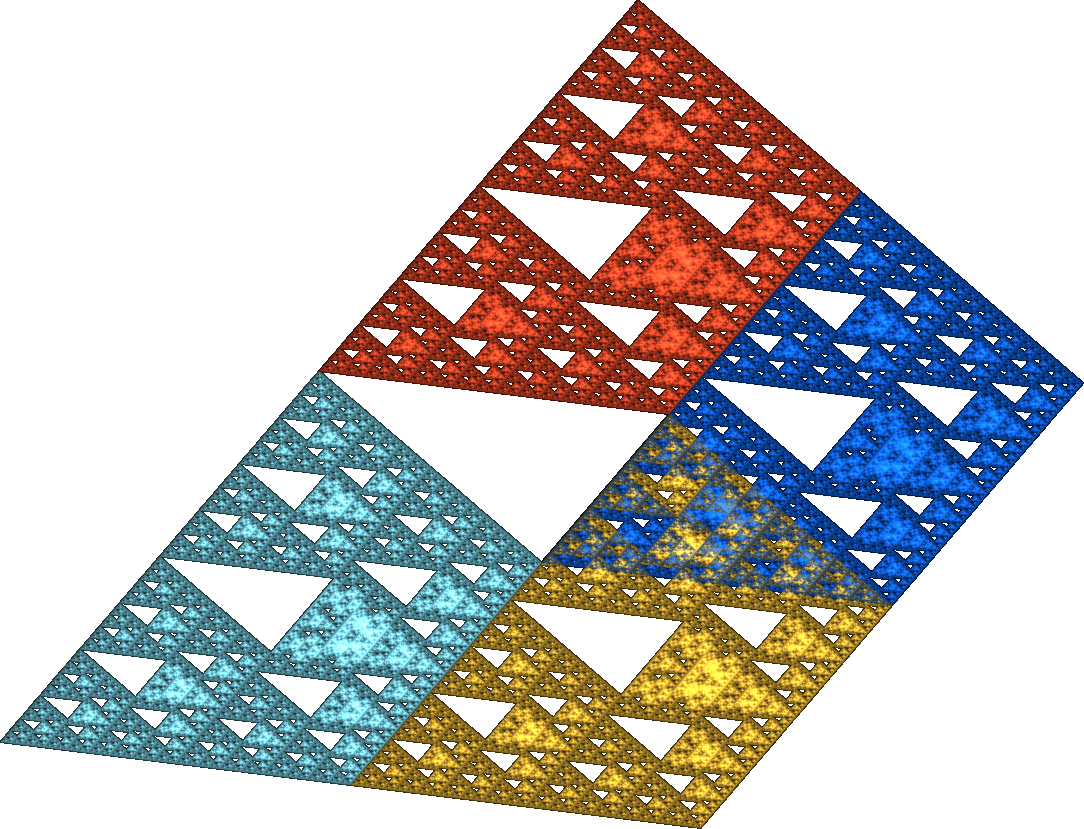}\quad
             \includegraphics[height=0.4\textwidth]{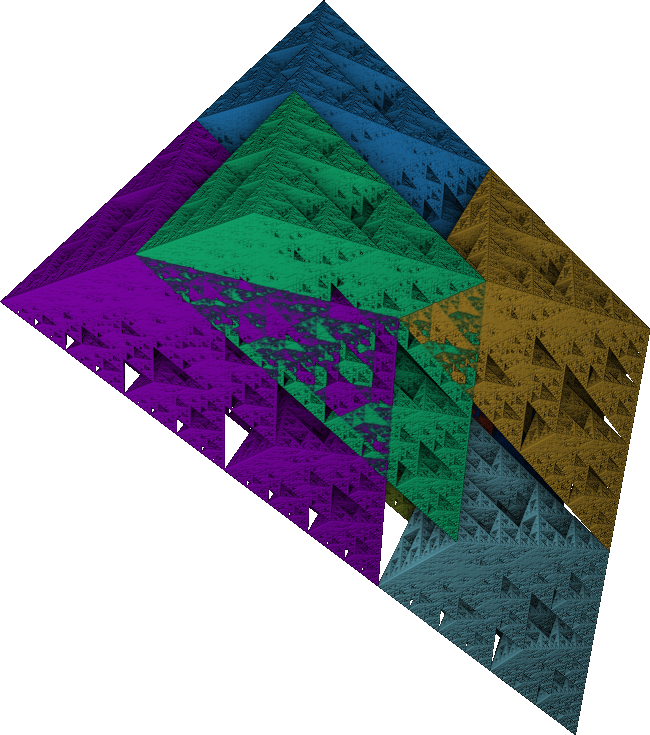}
	}
	\caption{Typical rep-tiles with complicated structure in plane and space.}
	\label{prob}
\end{figure}

\subsection*{Problems with the computer search.} Every kind of progress raises new problems. Our first experiments were disappointing. It can happen that 90 percent of the examples are cubes, which can be generated in many different data. This can be avoided by skipping examples with the same parameter values, cf. \cite{sierel}. The main problem, however, is that most of the generated examples have too complicated structure and bad geometric properties. Figure \ref{rep4} shows that plane reptiles can be disconnected: the rightmost example has two connected components. What is worse, the interior of this set has four components. The interior of the neighboring example has infinitely many components. Figure \ref{prob} shows another plane 4-rep-tile with fragmented interior, and a similar 8-rep-tile in space. Such sets fulfil the rep-tile definition but they cannot be physically realized as puzzle pieces.

Can we let the computer select the 25 most interesting examples in a list of 25 thousand?
Disconnected rep-tiles can be singled out by a simple algorithm. At present we have no method which controls the structure of the interior. However, some of the calculated parameters can be considered as measures of complexity. We may look for tiles which have boundary dimension 2 (that is, polygonal faces), and with a small (but not too small) number of neighbor types. They form a smaller collection which may be inspected by eyesight.

There are many options for the random search which will not be discussed here. We should keep in mind that the data space is huge. When the coordinates of the vector $v$ in \eqref{hk} are between -10 and +10, we have $n=48*21^3$ choices for each map. With eight maps there are $n^8\approx 10^{45}$ possible cases. Even many years' efforts could provide only a glimpse into our new territories.

\begin{figure}[h] 
	\centerline{\includegraphics[height=0.48\textwidth]{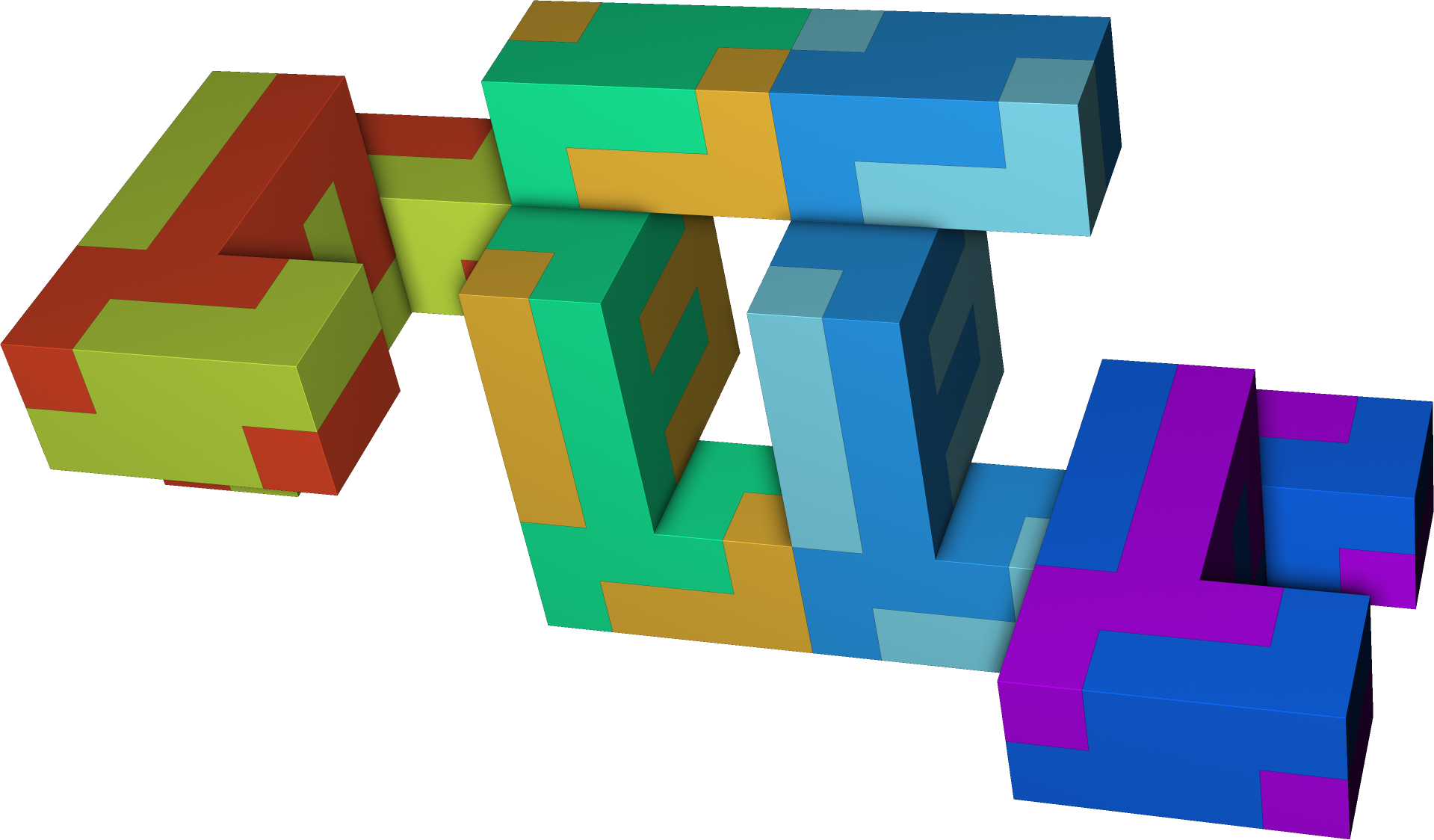}
	}
	\caption{Rep-tile with hole found by computer search.}
	\label{hole1}
\end{figure}

Perhaps we should be modest and study a smaller data space. For the cube family, we can consider $2\times 2\times 2$ rep-tiles as particular cases of 8-rep-tiles. We assume that the large set $g(A)$ is the union of two congruent blocks $C$ and $f_1(C),$ that $C$ is the union of two congruent blocks $D$ and $f_2(D)$, and that $D$ is the union of two copies of $A,$ which we call $f_3(A), f_4(A).$ The $f_k$ are again isometries of the form \eqref{hk}. Combining these equations, we obtain $g(A)$ as union of $f_3(A), f_4(A), f_2f_3(A), f_2f_4(A), f_1f_3(A), f_1f_4(A), f_1f_2f_3(A),$ and $f_1f_2f_4(A).$ This is a special case of equation \eqref{hut} which depends only on four maps. The corresponding data space includes $n^4\approx 4\cdot 10^{22}$ cases. It is still huge, but the chances not to get lost in our search are higher. Figure \ref{hole1} below was found with this approach.

\subsection*{Rep-tiles with holes.} Does there exist an $m$-rep-tile in space which is topologically equal to a torus? That is, its  interior is connected and has a single hole. The tile on the right of Figure \ref{rep8} has a kind of hole. But two of the little cubes which form the hole intersect only in an edge. So this is not really a solid ring: the interior of the tile has no hole. According to \cite{FS}, the question for tiles with hole was raised in 1998 by C. Goodman-Strauss and solved by G. van Ophuysen with an example for $m=24.$ A more general and abstract approach, with arbitrary number of holes and arbitrary large $m,$ was presented in \cite{coth}. We were interested in an example with $m=8$ and performed an extensive search of 8-rep-tiles. Among one million of examples, generated and pre-selected as described above, we found exactly one answer, shown in Figure \ref{hole1}. It is unlikely that anybody would find this example just by thinking and imagination!
The rep-tile can be made from wood, but it is mechanically impossible to assemble the pieces as shown in the figure.

Obviously, the tile consists of four congruent blocks. It is more difficult to see that the left and right part of the tile are also congruent, and that two small copies of the whole tile can be put together to form a block. The hole of the tile is realized by the two blocks in the middle. The blocks on the left and right are only needed to guarantee the self-similarity of the tile.

Can we solve the problem without the two superficial blocks? Consider only the middle part of Figure \ref{hole1}. When we consider a copy of this piece, rotated around 180 degrees, the copy will fill the hole, and both pieces together form a rectangular plate.  Moreover, the middle piece itself consists of four rectangular plates. We can vary their shapes in such a way that the previous sentence remains true. This leads to the following figure.  

\begin{figure}[h] 
	\centerline{\includegraphics[height=0.38\textwidth]{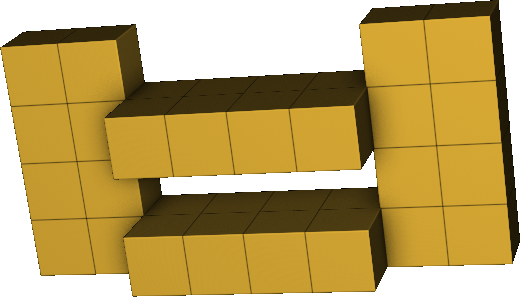}
	}
	\caption{Simplified rep-tile with hole.}
	\label{hole2}
\end{figure}

We show that this set is really an 8-rep-tile. It consists of four rectangular plates of size $4\times 2\times 1.$ Two copies of Figure \ref{hole2}, one of them rotated by 180 degrees, will fit together to form a rectangular plate of size $8\times 4\times 2.$ Thus eight copies of the set can be assembled to produce a similar set, expanded by the factor 2. This can be really done with material pieces, and the proof needs no calculation. It seems unlikely but not impossible that there is a still simpler rep-tile with a hole. 

What a happy end: man has shown to be stronger than machine, by finding a simpler tile. No:
both tiles are interesting, and  Figure \ref{hole1} was the starting point for Figure \ref{hole2}. Computers are here to stay, 
even in mathematical research. Let us use them - with critical interaction and new ideas.

\paragraph{Acknowledgment.}
The authors' cooperation was supported by the German Science Foundation (DFG), project Ba 1332/11-1.

\noindent
Christoph Bandt, 
Institute of Mathematics, 17487 Greifswald, Germany\\
\url{bandt@uni-greifswald.de}\\
Dmitry Mekhontsev,
Sobolev Institute of Mathematics, 630090 Novosibirsk, Russia\\
\url{mekhontsev@gmail.com}
\vfill\eject


\begin{thebibliography}{99}\small
\bibitem{Ba5} Christoph Bandt, Self-similar sets 5. {Integer} matrices and fractal tilings of
  {${\mathbb R}^n$}, \textit{Proc. Amer. Math. Soc.} 112 (1991), 549--562.
\bibitem{pinwheel} Christoph Bandt, Dmitry Mekhontsev and Andrei Tetenov, A single fractal pinwheel tile, 
\textit{Proc. Amer. Math. Soc.} 146 (2018), 1271--1285.
\bibitem{sierel} Christoph Bandt and Dmitry Mekhontsev, Elementary fractal geometry. New relatives of the Sierpi\'nski gasket, \textit{Chaos} \textbf{28} 063104  (2018).
\bibitem{Bar} Michael F. Barnsley, \textit{Fractals everywhere}, Academic Press, 2nd edition, 1993.
\bibitem{coth} Gregory R. Conner and J\"{o}rg M. Thuswaldner, Self-affine manifolds, \textit{Advances Math.} 289 (2016), 725--783.
\bibitem{FS} Dirk Frettloeh and Iwan Suschko, 3-Torus Rep-Tile, \url{http://www.eg-models.de/models/Polytopal_Complexes/2010.02.001/_direct_link.html}, 2010.
\bibitem{Ga} Martin Gardner, On rep-tiles, polygons that can make larger and smaller copies of themselves,
\textit{Scientific Amer.} 208 (1963) 154--164.
\bibitem{Gol} S.W. Golomb, Replicating figures in the plane, \textit{Math. Gaz.} 48 (1964) 403--412.
\bibitem{GM} Karl-Heinz Gr\"ochenig and W. Madych, Multiresolution analysis, Haar bases, and self-similar tilings,
\textit{IEEE Trans. Inform. Th.} 38 (2), Part 2 (1992) 558-568.
\bibitem{GS} Branko Gr\"{u}nbaum and G.C. Shephard, \textit{Patterns and Tilings},
Freeman, New York, 1987.
\bibitem{Ha} Herman Haverkort, No acute tetrahedron is an 8-reptile, arXiv:1508.03773v2 (2018)
\bibitem{Hu} John E. Hutchinson, Fractals and self-similarity,
\textit{Indiana University Mathematics Journal} 30 (1981) 713--747.
\bibitem{LJ} Anwei Liu and Barry Joe, On the shape of tetrahedra from bisection, 
\textit{Mathematics of Computation} 63 No. 207 (2013) 141-154.
\bibitem{MS11} Ji\v{r}i Matou\v{s}ek and Zuzana Safernov\'{a}, On the nonexistence of $k$-reptile tetrahedra,
\textit{Discrete Comput. Geom.} 46 (2011) 599--609.
\bibitem{M} Dmitry Mekhontsev, IFStile v1.7.4.4 (2018), \url{http://ifstile.com}
\bibitem{Se} Marjorie Senechal,  \textit{Quasicrystals and geometry}, Cambridge University Press, Cambridge 1995.
\bibitem{Vi} Andrew Vince, Rep-tiling Euclidean space,
\textit{Aequationes Math.} 50 (1995) 191--213.
\end{thebibliography}
\end{document}